\documentclass[twocolumn]{autart}    

\usepackage{graphicx}          

\usepackage{amsmath}
\usepackage{amssymb}

\newcommand*\diff{\mathop{}\!\mathrm{d}}%

\def\R{\mathbb{R}}

\usepackage{graphics} 
\usepackage{epsfig} 
\usepackage{epstopdf}
\usepackage{subfigure}
\usepackage{color}

\graphicspath{{./Figures/}}
\allowdisplaybreaks[4]

\newcounter{cst}

\begin{document}

\begin{frontmatter}

\title{Stability analysis of reaction-diffusion PDEs coupled at the boundaries with an ODE\thanksref{footnoteinfo}} 

\thanks[footnoteinfo]{Corresponding author H.~Lhachemi. The work of C. Prieur has been partially supported by MIAI@Grenoble Alpes (ANR-19-P3IA-0003)
}

\author[CS]{Hugo Lhachemi}\ead{hugo.lhachemi@centralesupelec.fr},
\author[GIPSA-lab]{Christophe Prieur}\ead{christophe.prieur@gipsa-lab.fr}, 

\address[CS]{Universit{\'e} Paris-Saclay, CNRS, CentraleSup{\'e}lec, Laboratoire des signaux et syst{\`e}mes, 91190, Gif-sur-Yvette, France}  
\address[GIPSA-lab]{Universit{\'e} Grenoble Alpes, CNRS, Grenoble-INP, GIPSA-lab, F-38000, Grenoble, France}             

\begin{keyword}                           
Coupled PDE-ODE, stability, reaction-diffusion equation, modal decomposition, LMI              
\end{keyword}                             

\begin{abstract}                          
This paper addresses the derivation of generic and tractable sufficient conditions ensuring the stability of a coupled system composed of a reaction-diffusion partial differential equation (PDE) and a finite-dimensional linear time invariant ordinary differential equation (ODE). The coupling of the PDE with the ODE is located either at the boundaries or in the domain of the reaction-diffusion equation and takes the form of the input and output of the ODE. We investigate boundary Dirichlet/Neumann/Robin couplings, as well as in-domain Dirichlet/Neumann couplings. The adopted approach relies on the spectral reduction of the problem by projecting the trajectory of the PDE into a Hilbert basis composed of the eigenvectors of the underlying Sturm-Liouville operator and yields a set of sufficient stability conditions taking the form of LMIs. We propose numerical examples, consisting of an unstable reaction-diffusion equation and an unstable ODE, such that the application of the derived stability conditions ensure the stability of the resulting coupled PDE-ODE system. 
\end{abstract}

\end{frontmatter}

\section{Introduction}

The stability analysis and control of coupled PDE-ODE systems has emerged relatively recently in the literature (and more generaly PDEs with dynamical boundary conditions, see e.g.~\cite{nicaise2009stability}). Such a trend was driven by a certain number of practical applications involving a finite-dimensional dynamics coupled with a phenomenon described by a PDE. This includes, to cite a few, solid–gas interaction of heat diffusion and chemical reaction~\cite{tang2011state}, flexible cranes~\cite{he2016cooperative}, flexible aircraft~\cite{lhachemi2018boundary}, drilling mechanisms~\cite{barreau2018exponential}, and power converters connected to transmission lines~\cite{daafouz2014nonlinear}. PDE-ODE coupling can also arise due to feedback control. Indeed, the PDE can represent the open-loop plant to be controlled while the ODE part gathers controller and actuator dynamics, see e.g.~\cite{krstic2009compensating,susto2010control}. Conversely, the PDE can represent the dynamics of an actuator (e.g., heat or flux sensors) that is embedded into the closed-loop control of a finite-dimensional plant modeled by an ODE.

The stabilization of PDE-ODE couplings has attracted much attention in the recent years. One of the very first contributions in this field was reported in~\cite{krstic2009compensating} dealing with the state-feedback stabilization and the observer design of a diffusion PDE cascaded with an ODE via Dirichlet connection (see also~\cite{susto2010control} for the case of Neumann interconnections). Such a problem can be interpreted as a compensation problem of an infinite-dimensional input dynamics~\cite{krstic2010compensation} and was solved by employing a backstepping control design procedure. This approach was also reported in the case of string equation in~\cite{krstic2009compensatingString,susto2010control} and was later on applied to other types of PDEs such as beam~\cite{wu2014static} and linearized Korteweg–de Vries~\cite{ayadi2018exponential} equations. This backstepping-based procedure for PDE-ODE cascades was then extended to other boundary stabilization problems such as wave PDEs cascaded with MIMO LTI systems~\cite{bekiaris2010compensating}, a diffusion PDE coupled with an ODE~\cite{tang2011state}, and a diffusion PDE sandwiched between two ODEs~\cite{wang2019output}. The robustness of certain of these control strategies for the stabilization of PDE-ODE cascades were studied in~\cite{krstic2009compensating,susto2010control,sanz2018robust}, particularly for heat equations w.r.t. the diffusion coefficient and the length of the diffusion domain. Other extensions embracing the augmentation of the backstepping transformation with either adaptive or sliding mode control have been investigated in~\cite{li2014adaptive,wang2015sliding}. Recently, a different control design strategy using Sylvester equation was proposed in~\cite{natarajan2021compensating} for ODE-PDE and PDE-ODE cascades.

In this context, the focus of this paper is put on the derivation of generic and tractable sufficient stability conditions ensuring the stability of coupled systems composed of a general reaction-diffusion PDE and an ODE, with various coupling configurations, rather than on the design of a particular controller for a specific setting. The derivation of such analysis tools for stability assessment of PDE-ODE loops is of primary importance. Indeed, such generic stability conditions can be used to assess \emph{a posteriori} the stability of open-loop unstable reaction-diffusion PDEs when placed in closed-loop with a controller designed empirically using a finite-dimensional truncated model of the PDE. Conversely, considering the infinite-dimensional control strategies reported in the previous paragraph, their practical implementation require their finite-dimensional approximation. In that case, stability analysis tools are required to assess that the finite-dimensional approximation of the controller dynamics still achieves the stabilization of the PDE.

The traditional approach for studying the stability of coupled PDE-ODE systems consists of the adequate selection of a Lyapunov functional. At a very high level, the general trend is to build the Lyapunov functional by considering terms related to 1) the energy of the PDE (measured via a relevant norm); 2) the energy of the ODE; 3) the coupling of the PDE-ODE system. Such Lyapunov functionals can be built manually~\cite{coron2004global,coron2006global,lhachemi2020feedback} but can also be obtained numerically by considering very general Lyapunov functional candidates while resorting to numerical methods, such as a sum of square procedure, to obtain an admissible suitable set of parameters~\cite{ahmadi2016dissipation,gahlawat2016convex}.

In the abovementioned context, a number of contributions have been reported in the recent years to study the stability of coupled PDE-ODE systems with couplings occurring at the boundaries of the PDE. A first fruitful approach relies on the introduction of a partial integral representation of the PDE~\cite{peet2021partial} in order to study the stability of PDE-ODE loops using linear matrix inequalities (LMIs). Such an approach can be used to study PDE-ODE loops using convex optimization tools~\cite{das2020robust,peet2021representation,shivakumar2019generalized}. A second fruitful approach relies on the use of Legendre polynomials as a basis of projection for the PDE trajectories. In essence, this consists of the construction of a classical Lyapunov functional accounting for the PDE and ODE parts considered separately while adding a cross quadratic term mixing the state of the ODE with a finite number of coefficients of projection of the PDE trajectory into the basis of Legendre polynomials. Such an approach was reported in~\cite{baudouin2019stability} for the study of a coupled system composed of a reaction PDE and an ODE. This method was also reported in~\cite{barreau2018lyapunov} in the case of a string equation coupled with an ODE, as well as in~\cite{barreau2018input} for the study of input-output stability. Input-output stability properties for coupled PDE-ODE systems using Legendre polynomials-based projections was further investigated in~\cite{matthieu2020integral}. Finally, the stability of abstract boundary control systems with dynamic boundary conditions and positive underlying $C_0$-semigroups was studied in~\cite{boulouz2021well}. 

In this paper, we study the stability of a generic 1-D reaction diffusion equation coupled with a finite-dimensional ODE. The approach adopted in this work differs from the methods described in the previous paragraph because it relies on spectral reduction methods. These spectral reduction methods are used to build a suitable Lyapunov functional candidate and derive a set of tractable LMI conditions ensuring the exponential stability of the coupled PDE-ODE system. Compared to~\cite{baudouin2019stability}, which was concerned with an open-loop stable constant coefficient diffusion PDE with left and right Dirichlet couplings, our approach allows the consideration of generic reaction-diffusion PDEs that are possibly open-loop unstable and with variety of couplings that include Dirichlet, Neumann, and Robin traces. Compared to~\cite{das2020robust,peet2021representation,shivakumar2019generalized} the approach adopted in this paper allows the coupling of the ODE with the PDE through a Dirichlet/Neumann trace that can be located either at the boundary or inside the spatial domain. Moreover, the exponential stability results derived in this paper are established for system trajectories evaluated in $H^1$-norm. This feature has two important implications: 1) the exponential decrease of the PDE trajectories in $L^\infty$-norm\footnote{This result immediately follows from our stability result established in $H^1$-norm and the fact that the $L^\infty$-norm is bounded by the $H^1$-norm. Note however that this does not imply stability in $L^\infty$-norm.}; and 2) the exponential decay of the coupling channels between ODE and PDE components. This last point is of paramount importance for practical applications because it ensures that the signals in the actuation/sensing channels are also convergent. The relevance of these LMI conditions are assessed based on numerical examples associated with PDEs and ODEs that are all unstable.

The rest of the paper is organized as follows. Section~\ref{sec: preliminaries} describes the notations and reports a number of basic properties for Sturm-Liouville operators. Then the study is split into two parts. Firstly, the case of a Dirichlet trace used as an input for the ODE is investigated in Section~\ref{sec: Dirichlet trace as an input of the ODE}. Secondly, the case of a Neumann trace used as an input for the ODE is reported in Section~\ref{sec: Neumann trace as an input of the ODE}. Finally, concluding remarks are formulated in Section~\ref{sec: conclusion}.

\section{Notation and properties}\label{sec: preliminaries}

Spaces $\R^n$ are endowed with the Euclidean norm denoted by $\Vert\cdot\Vert$. The associated induced norms of matrices are also denoted by $\Vert\cdot\Vert$. $L^2(0,1)$ stands for the space of square integrable functions on $(0,1)$ and is endowed with the inner product $\langle f , g \rangle = \int_0^1 f(x) g(x) \,\mathrm{d}x$ and the norm is denoted by $\Vert \cdot \Vert_{L^2}$. For an integer $m \geq 1$, the $m$-order Sobolev space is denoted by $H^m(0,1)$ and is endowed with its usual norm $\Vert \cdot \Vert_{H^m}$. For a symmetric matrix $P \in\R^{n \times n}$, $P \succeq 0$ (resp. $P \succ 0$) means that $P$ is positive semi-definite (resp. positive definite) while $\lambda_M(P)$ (resp. $\lambda_m(P)$) denotes its maximal (resp. minimal) eigenvalue.

Let $\theta_1,\theta_2 \in [0,\pi/2]$, $p \in \mathcal{C}^1([0,1])$ with $p > 0$, and $q \in \mathcal{C}^0([0,1])$ with $q \geq 0$. Let the Sturm-Liouville operator $\mathcal{A} : D(\mathcal{A}) \subset L^2(0,1) \rightarrow L^2(0,1)$ be defined by $\mathcal{A} f = - (pf')' + q f$ on the domain $D(\mathcal{A}) = \{ f \in H^2(0,1) \,:\, \cos(\theta_1) f(0) - \sin(\theta_1) f'(0) = 0 ,\, \allowbreak \cos(\theta_2) f(1) + \sin(\theta_2) f'(1) = 0 \}$. The operator $\mathcal{A}$ is self-adjoint and its eigenvalues $\lambda_n$, $n \geq 1$, are simple, non negative, and form an increasing sequence with $\lambda_n \rightarrow + \infty$ as $n \rightarrow + \infty$. Moreover, the associated unit eigenvectors $\phi_n \in L^2(0,1)$ form a Hilbert basis and we also have $D(\mathcal{A}) = \{ f \in L^2(0,1) \,:\, \sum_{n\geq 1} \vert \lambda_n \vert ^2 \vert \left< f , \phi_n \right> \vert^2 < +\infty \}$ with $\mathcal{A}f = \sum_{n \geq 1} \lambda_n \left< f, \phi_n \right> \phi_n$.

Let $p_*,p^*,q^* \in \R$ be such that $0 < p_* \leq p(\xi) \leq p^*$ and $0 \leq q(\xi) \leq q^*$ for all $\xi \in [0,1]$, then it holds~\cite{orlov2017general}:
\begin{equation}\label{eq: estimation lambda_n}
0 \leq \pi^2 (n-1)^2 p_* \leq \lambda_n \leq \pi^2 n^2 p^* + q^*
\end{equation}
for all $n \geq 1$. Assuming further than $p \in \mathcal{C}^2([0,1])$, we have for any given $\xi \in [0,1]$ that $\phi_n (\xi) = O(1)$ and $\phi_n'(\xi) = O(\sqrt{\lambda_n})$ as $n \rightarrow + \infty$ (see~\cite{orlov2017general}), where $O$ denotes the Bachmann–Landau asymptotic notation. Moreover, we also have, for all $f \in D(\mathcal{A})$, $\left< \mathcal{A} f , f \right> = \sum_{n \geq 1} \lambda_n \left< f , \phi_n \right>^2 = \int_0^1 p (f')^2 + q f^2 \,\mathrm{d}\xi + p(0) f(0) f'(0) - p(1) f(1) f'(1)$. Hence, provided $q > 0$ and because $\theta_1,\theta_2 \in [0,\pi/2]$ and $p > 0$, we obtain the existence of constants $C_1,C_2 > 0$ so that
\begin{align}
C_1 \Vert f \Vert_{H^1}^2 \leq 
\sum_{n \geq 1} \lambda_n \left< f , \phi_n \right>^2
= \left< \mathcal{A}f , f \right>
\leq C_2 \Vert f \Vert_{H^1}^2 \label{eq: inner product Af and f}
\end{align}
for any $f \in D(\mathcal{A})$. This in particular implies that, for any $f \in D(\mathcal{A})$ and any $\xi \in [0,1]$, $f(\xi) = \sum_{n \geq 1} \left< f , \phi_n \right> \phi_n(\xi)$ and $f'(\xi) = \sum_{n \geq 1} \left< f , \phi_n \right> \phi_n'(\xi)$. For any $\alpha \in (0,1)$, we introduce the fractional powers of $\mathcal{A}$ by defining $D(\mathcal{A}^\alpha) = \{ f \in L^2(0,1) \,:\, \sum_{n\geq 1} \vert \lambda_n \vert ^{2\alpha} \vert \left< f , \phi_n \right> \vert^2 < +\infty \}$ and $\mathcal{A}^{\alpha} f = \sum_{n \geq 1} \lambda_n^{\alpha} \left< f , \phi_n \right> \phi_n$. 

For any $f \in L^2(0,1)$ we define $\mathcal{P}_N f = f - \sum_{n=1}^{N} \left<f,\phi_n\right>\phi_n = \sum_{n\geq N+1} \left<f,\phi_n\right>\phi_n$

\section{Dirichlet trace as an input of the ODE}\label{sec: Dirichlet trace as an input of the ODE}

\subsection{Coupled PDE-ODE systems}

We consider in this section the following PDE-ODE system:
\begin{subequations}\label{eq: coupled PDE-ODE system 1}
\begin{align}
& z_t(t,\xi) = (p z_\xi)_\xi(t,\xi) - \tilde{q}(\xi) z(t,\xi) \\
& \cos(\theta_1) z(t,0) - \sin(\theta_1) z_\xi(t,0) = 0 \label{eq: coupled PDE-ODE system 1 - BC1} \\
& \cos(\theta_2) z(t,1) + \sin(\theta_2) z_\xi(t,1) = y(t) = C x(t) \label{eq: coupled PDE-ODE system 1 - BC2} \\
& \dot{x}(t) = A x(t) + Bz(t,\zeta_m) \label{eq: coupled PDE-ODE system 1 - ODE} \\
& z(0,\xi) = z_0(\xi) , \quad x(0) = x_0
\end{align}
\end{subequations}
for $t > 0$ and $\xi\in(0,1)$ where $\theta_1,\theta_2 \in [0,\pi/2]$, $p \in \mathcal{C}^2([0,1])$ with $p > 0$, $\tilde{q} \in \mathcal{C}^0([0,1])$, and $\zeta_m \in[0,1]$. Here $A\in\R^{n \times n}$, $B\in\R^{n}$, and $C\in\R^{1 \times n}$ are matrices, $z_0 \in L^2(0,1)$ and $x_0 \in\R^n$ are initial conditions, and $z(t,\cdot)\in L^2(0,1)$ and $x(t) \in\R^n$ are the state of the reaction-diffusion PDE and of the ODE at time $t$, respectively. 

The PDE-ODE system (\ref{eq: coupled PDE-ODE system 1}) consists of a reaction-diffusion PDE coupled with an ODE. The output $y(t) = Cx(t)$ of the ODE is seen as a boundary input for the PDE and is applied at the right Robin boundary condition. Conversely, the pointwise Dirichlet trace $z(t,\zeta_m)$ is seen as an input of the ODE (\ref{eq: coupled PDE-ODE system 1 - ODE}). The objective of this section is to derive numerically tractable sufficient conditions ensuring the exponential stability of the PDE-ODE system (\ref{eq: coupled PDE-ODE system 1}) when evaluating the PDE trajectory in $H^1$-norm.

We introduce without loss of generality $q \in \mathcal{C}^0([0,1])$ and $q_c \in \R$ such that 
\begin{equation}\label{eq: decomposition of the reaction term}
\tilde{q} = q - q_c , \quad q > 0.
\end{equation}

\begin{rem}
Even if the presentation focuses on the case $\theta_1,\theta_2 \in [0,\pi/2]$, the derived results can be extended to $\theta_1,\theta_2 \in [0,2\pi]$. Indeed, considering first the case $\theta_1,\theta_2 \in [0,\pi]$, the proposed strategy also applies provided 1) $q > 0$ from (\ref{eq: decomposition of the reaction term}) is selected large enough so that the estimates (\ref{eq: inner product Af and f}) still hold for some constants $C_1,C_2 > 0$; 2) the change of variable (\ref{eq: change of variable 1}) is replaced by $w(t,\xi) = z(t,\xi) - \frac{\xi^\alpha}{\cos\theta_2 + \alpha \sin\theta_2} y(t)$ for any fixed $\alpha > 1$ selected so that $\cos(\theta_2) + \alpha \cos(\theta_2) \neq 0$. Finally, in view of (\ref{eq: coupled PDE-ODE system 1 - BC1}-\ref{eq: coupled PDE-ODE system 1 - BC2}), the general case $\theta_1,\theta_2 \in [0,2\pi]$ reduces to the case $\theta_1,\theta_2 \in [0,\pi]$ by proceeding with the following substitutions: 1) if $\theta_1 \in (\pi,2\pi]$ then $\theta_1 \leftarrow \theta_1 - \pi$; 2) if $\theta_2 \in (\pi,2\pi]$ then $\theta_2 \leftarrow \theta_2 - \pi$ and $C \leftarrow -C$.
\end{rem}

\begin{rem}
System (\ref{eq: coupled PDE-ODE system 1}), as well as system (\ref{eq: coupled PDE-ODE system 3}) that will be described in the next section, can be used to represent a variety of practical situations. For instance the PDE part can stand for a reaction-diffusion process coupled with a finite-dimensional LTI controller materialized by the ODE. Conversely, the ODE part can merge both finite-dimensional LTI plant along with its associated finite-dimensional LTI controller while the PDE part describes the sensor dynamics. This latter situation is similar to the one described in~\cite{krstic2009compensating} where a controller was designed for a cascaded PDE-ODE system using a backstepping transformation. One of the main motivations for deriving generic stability conditions for coupled PDE-ODE systems such as (\ref{eq: coupled PDE-ODE system 1}) is ultimately when the to-be-implemented finite-dimensional controller is computed either on a finite-dimensional approximation of the PDE or via the approximation of an infinite-dimensional output feedback controller (obtained, e.g., using backstepping control design procedures).
\end{rem}

\begin{rem}
The PDE-ODE system (\ref{eq: coupled PDE-ODE system 1}) with $\theta_1 = \pi/2$, $\theta_2 = 0$, and $\zeta_m = 0$ was studied in~\cite{baudouin2019stability} in the case of a stable diffusion PDE (i.e., without reaction term) and with a constant diffusion coefficient using the projection of the PDE trajectories into a finite subset of Legendre polynomials.
\end{rem}

\subsection{Preliminary spectral reduction}

We rewrite (\ref{eq: coupled PDE-ODE system 1}) under an equivalent PDE-ODE system with homogeneous boundary conditions. Specifically, introducing the change of variable
\begin{equation}\label{eq: change of variable 1}
w(t,\xi) = z(t,\xi) - \frac{\xi^2}{\cos\theta_2 + 2\sin\theta_2} y(t) 
\end{equation}
we infer that (\ref{eq: coupled PDE-ODE system 1}) is equivalent to
\begin{subequations}\label{eq: homogeneous coupled PDE-ODE system 1}
\begin{align}
& w_t(t,\xi) = (p w_\xi)_\xi(t,\xi) + (q_c-q(\xi)) w(t,\xi) \\
& \phantom{w_t(t,\xi) =}\; + a(\xi) y(t) + b(\xi) \dot{y}(t) \nonumber \\
& \cos(\theta_1) w(t,0) - \sin(\theta_1) w_\xi(t,0) = 0 \\
& \cos(\theta_2) w(t,1) + \sin(\theta_2) w_\xi(t,1) = 0 \\
& y(t)=Cx(t) \\
& \dot{x}(t) = A x(t) + B \left( w(t,\zeta_m \right) + \mu_m y(t) )  \\
& w(0,\xi) = w_0(\xi) , \quad x(0) = x_0
\end{align}
\end{subequations}
with $a(\xi) = \frac{1}{\cos\theta_2+2\sin\theta_2} \{ 2p(\xi) + 2 \xi p'(\xi) + (q_c-q(\xi)) \xi^2 \}$, $b(\xi) = -\frac{\xi^2}{\cos\theta_2+2\sin\theta_2}$, $\mu_m = - b(\zeta_m)$, and $w_0(\xi) = z_0(\xi) - \frac{\xi^2}{\cos\theta_2+2\sin\theta_2} y(0)$. 

After this change of variable, the well-posedness in terms of classical solutions of the above PDE-ODE system for initial conditions $w_0 \in D(\mathcal{A}^{1/2})$ and $x_0 \in\R^n$ is a consequence of~\cite[Thm.~6.3.1 and~6.3.3]{pazy2012semigroups}. More precisely, we have for any $w_0 \in D(\mathcal{A}^{1/2})$ and any $x_0 \in\R^n$ the existence and uniqueness of a classical solution $(w,x) \in \mathcal{C}^0([0,\infty);L^2(0,1) \times \R^n) \cap \mathcal{C}^1((0,\infty);L^2(0,1) \times \R^n)$ with $w(t,\cdot) \in D(\mathcal{A})$ for all $t > 0$. Moreover, from the proof of ~\cite[Thm.~6.3.1]{pazy2012semigroups}, we have $\mathcal{A} w \in \mathcal{C}^0((0,\infty);L^2(0,1))$ and $\mathcal{A}^{1/2}w \in \mathcal{C}^0([0,\infty);L^2(0,1))$.

We now introduce the Hilbert basis $\{\phi_i \,:\, i \geq 1\}$ of $L^2(0,1)$ formed by the eigenvectors of the Sturm-Liouville operator $\mathcal{A}$. We introduce the coefficients of projection:
\begin{equation}\label{eq: coefficients of projection}
w_i(t) = \left<w(t,\cdot),\phi_i\right> , \quad
a_i = \left<a,\phi_i\right> , \quad
b_i = \left<b,\phi_i\right> 
\end{equation}
and define $c_i = \phi_i(\zeta_m)$ for $i \geq 1$. Considering classical solutions, we obtain that 
\begin{subequations}\label{eq: spectral reduction 1}
\begin{align}
& \dot{w}_i(t) = (-\lambda_i + q_c) w_i(t) + a_i C x(t) \label{eq: spectral reduction 1 - 1} \\
& \phantom{\dot{w}_i(t) =}\; + b_i C \left\{ ( A + \mu_m BC ) x(t) + B \sum_{j \geq 1} c_j w_j(t) \right\} \nonumber \\
& \dot{x}(t) = ( A + \mu_m BC ) x(t) + B \sum_{j \geq 1} c_j w_j(t) \label{eq: spectral reduction 1 - 2}
\end{align}
\end{subequations}
for $i \geq 1$. The adopted stability analysis procedure relies now on the introduction of a finite dimensional model that captures the dynamics of the ODE (\ref{eq: spectral reduction 1 - 2}) along with the $N \geq 1$ first modes $w_i$ of the PDE plant, described by (\ref{eq: spectral reduction 1 - 1}), while bounding the effect of the residue of measurement $\mathcal{R}(t) = \sum_{i \geq N+1} c_i w_i(t)$ by using Lyapunov's direct method. To do so, we define 
\begin{align*}
W(t) & = \begin{bmatrix} w_1(t) & \ldots & w_N(t) \end{bmatrix}^\top \in\R^N , \\
A_N & = \mathrm{diag}(-\lambda_1+q_c,\ldots,-\lambda_N+q_c) \in\R^{N \times N} , \\
B_{a,N} & = \begin{bmatrix} a_1 & \ldots & a_N \end{bmatrix}^\top \in\R^N , \\
B_{b,N} & = \begin{bmatrix} b_1 & \ldots & b_N \end{bmatrix}^\top \in\R^N , \\
C_N & = \begin{bmatrix} c_1 & \ldots & c_N \end{bmatrix} \in\R^{1 \times N} .
\end{align*}
We infer from (\ref{eq: spectral reduction 1 - 1}) that 
\begin{align*}
\dot{W}(t) 
& = ( A_N + B_{b,N} C B C_N ) W(t)  + B_{b,N} C B \mathcal{R}(t) \\
& \phantom{=}\; + ( B_{a,N} C + B_{b,N} C ( A + \mu_m BC ) ) x(t) .
\end{align*}
Combining this latter identity with (\ref{eq: spectral reduction 1 - 2}) while defining
\begin{equation*}
X(t) = \begin{bmatrix} W(t) \\ x(t) \end{bmatrix} \in\R^{N+n},
\end{equation*}
we infer that
\begin{equation}\label{eq: truncated model 1}
\dot{X}(t) = F X(t) + G \mathcal{R}(t)
\end{equation}
where
\begin{equation*}
F = \begin{bmatrix} A_N + B_{b,N}CBC_N & B_{a,N}C+B_{b,N}C( A + \mu_m BC ) \\ BC_N & A + \mu_m BC \end{bmatrix}
\end{equation*}
and
\begin{equation*}
G = \begin{bmatrix} B_{b,N}CB \\ B \end{bmatrix} .
\end{equation*}
Hence, the ODE (\ref{eq: truncated model 1}) describes the dynamics of the ODE and of the $N$ first modes $w_i$ of the PDE plant while taking as an input the residue of measurement $\mathcal{R}(t) = \sum_{i \geq N+1} c_i w_i(t)$. The residual dynamics, which corresponds to the modes $i \geq N+1$, is characterized by
\begin{align}
\dot{w}_i(t)& = (-\lambda_i + q_c) w_i(t) + a_i C x(t) + b_i C ( A + \mu_m BC ) x(t) \nonumber \\
& \phantom{=}\; + b_i CBC_N W(t) + b_i CB \mathcal{R}(t) . \label{eq: residual dynamics 1}  
\end{align}

In preparation of the stability analysis, we introduce the matrix 
\begin{equation*}
H = \begin{bmatrix} H_{1,1} & 0 \\ 0 & H_{2,2} \end{bmatrix}
\end{equation*}
with $H_{1,1} = \Vert \mathcal{P}_N b \Vert_{L^2}^2 C_N^\top B^\top C^\top C B C_N$ and $H_{2,2} = \Vert \mathcal{P}_N a \Vert_{L^2}^2 C^\top C + \Vert \mathcal{P}_N b \Vert_{L^2}^2 ( A + \mu_m BC )^\top C^\top C ( A + \mu_m BC )$. We finally define the constant defined by $M_{1,\phi} = \sum_{i \geq N+1} \frac{\phi_i(\zeta_m)^2}{\lambda_i}$ which is finite because (\ref{eq: estimation lambda_n}) along with $\phi_n (\zeta_m) = O(1)$ as $n \rightarrow + \infty$.

\subsection{Main result}

We can now introduce the main result of this section.

\begin{thm}\label{thm1}
Let $\theta_1,\theta_2 \in [0,\pi/2]$, $p \in \mathcal{C}^2([0,1])$ with $p > 0$, $\tilde{q} \in \mathcal{C}^0([0,1])$, $\zeta_m \in[0,1]$, $A\in\R^{n \times n}$, $B\in\R^{n}$, and $C\in\R^{1 \times n}$ be given. Let $q \in \mathcal{C}^0([0,1])$ and $q_c \in\R$ be such that (\ref{eq: decomposition of the reaction term}) holds. Assume that there exist $N \geq 1$, $P \succ 0$, $\alpha > 2$, and $\beta > 0$ such that $\Theta_1 \prec 0$ and $\Theta_2 \prec 0$ where
\begin{align*}
\Theta_1 & = 
\begin{bmatrix} F^\top P + P F + \alpha H & P G \\ G^\top P & \alpha \Vert \mathcal{P}_N b \Vert_{L^2}^2 (CB)^2 - \beta \end{bmatrix} , \\
\Theta_2 & = 
\begin{bmatrix} -\lambda_{N+1}+q_c+\frac{\beta M_{1,\phi}}{2} & \sqrt{2\lambda_{N+1}} \\ \sqrt{2\lambda_{N+1}} & -\alpha \end{bmatrix} .
\end{align*} 
Then there exist constants $\eta,M > 0$ such that, for any initial conditions $z_0 \in H^2(0,1)$ and $x_0 \in \R^n$ such that $\cos(\theta_1) z_0(0) - \sin(\theta_1) z_0'(0) = 0$ and $\cos(\theta_2) z_0(1) + \sin(\theta_2) z_0'(1)=Cx_0$, the classical solution of (\ref{eq: coupled PDE-ODE system 1}) satisfies 
\begin{subequations}\label{eq: stab property 1}
\begin{equation}
\Vert z(t,\cdot) \Vert_{H^1}^2 + \Vert x(t) \Vert^2 \leq M e^{-2\eta t} ( \Vert z_0 \Vert_{H^1}^2 + \Vert x_0 \Vert^2 )
\end{equation}
with coupling channels such that
\begin{equation}
z(t,\zeta_m)^2 + y(t)^2 \leq  M e^{-2\eta t} ( \Vert z_0 \Vert_{H^1}^2 + \Vert x_0 \Vert^2 )
\end{equation}  
\end{subequations}
for all $t > 0$.
\end{thm}

\begin{rem}
The conclusions of Theorem~\ref{thm1} actually hold for any initial conditions such that $w_0 \in D(\mathcal{A}^{1/2})$. For example the case $\theta_1 = \pi/2$ and $\theta_2 = 0$ leads to $D(\mathcal{A}^{1/2}) = \left\{ f \in H^1(0,1) \,:\, f(1)=0 \right\}$. In this setting the conclusions of Theorem~\ref{thm1} hold for any $z_0 \in H^1(0,1)$ and $x_0 \in \R^n$ such that $z_0(1) = Cx_0$. Similarly, the case $\theta_1 = \theta_2 = \pi/2$ gives $D(\mathcal{A}^{1/2}) = H^1(0,1)$, implying that the conclusions of Theorem~\ref{thm1} hold for any $z_0 \in H^1(0,1)$ and $x_0 \in \R^n$.
\end{rem}

\textbf{Proof.}
Let $N \geq 1$, $P \succ 0$, $\alpha > 2$, and $\beta > 0$ such that $\Theta_1 \prec 0$ and $\Theta_2 \prec 0$. Hence, there exist $\eta > 0$ such that $\Theta_{1,\eta} \preceq 0$ and $\Theta_{2,\eta} \preceq 0$ where
\begin{align*}
\Theta_{1,\eta} & = 
\begin{bmatrix} F^\top P + P F + 2 \eta P + \alpha H & P G \\ G^\top P & \alpha \Vert \mathcal{P}_N b \Vert_{L^2}^2 (CB)^2 - \beta \end{bmatrix} , \\
\Theta_{2,\eta} & = 
\begin{bmatrix} -\lambda_{N+1}+q_c+\eta+\frac{\beta M_{1,\phi}}{2} & \sqrt{2\lambda_{N+1}} \\ \sqrt{2\lambda_{N+1}} & -\alpha \end{bmatrix} .
\end{align*} 
Define the Lyapunov functionnal candidate
\begin{equation}\label{eq: Lyapunov functionnal}
V(X,w) = X^\top P X + \sum_{i \geq N+1} \lambda_i \left< w , \phi_i \right>^2
\end{equation}
with $X\in\R^{N+n}$ and $w \in D(\mathcal{A})$. The first term of the above functional accounts for the finite-dimensional truncated model (\ref{eq: truncated model 1}) while the series is used to study the stability of the residual dynamics described by (\ref{eq: residual dynamics 1}) and to bound the effect of the residue of measurement $\mathcal{R}(t) = \sum_{i \geq N+1} c_i w_i(t)$ which is acting as an input of (\ref{eq: truncated model 1}). With the slight abuse of notation $V(t) = V(X(t),w(t))$, the computation of the time derivative of $V$ along the system trajectories (\ref{eq: truncated model 1}) and (\ref{eq: residual dynamics 1}) gives for $t > 0$
\begin{align*}
& \dot{V}(t) = 2 X(t)^\top P \dot{X}(t) + 2 \sum_{i \geq N+1} \lambda_i w_i(t) \dot{w}_i(t) \\
& = X(t)^\top ( F^\top P + P F ) X(t) + 2 X(t)^\top P G \mathcal{R}(t)  \\
& \phantom{=}\; + 2 \sum_{i \geq N+1} \lambda_i (-\lambda_i+q_c) w_i(t)^2 \\
& \phantom{=}\; + 2 \sum_{i \geq N+1} \lambda_i w_i(t) a_i C x(t) \\
& \phantom{=}\;  + 2 \sum_{i \geq N+1} \lambda_i w_i(t) b_i C ( A + \mu_m BC ) x(t) \\
& \phantom{=}\; + 2 \sum_{i \geq N+1} \lambda_i w_i(t) b_i C B C_N W(t) \\
& \phantom{=}\; + 2 \sum_{i \geq N+1} \lambda_i w_i(t) b_i C B \mathcal{R}(t) .
\end{align*}
We estimate the four latter series by using Young's inequality. For instance, the first term is estimated as
\begin{align*}
& 2 \sum_{i \geq N+1} \lambda_i w_i(t) a_i C x(t) \\
& \leq \sum_{i \geq N+1} \left\{ \frac{1}{\alpha}\lambda_{i}^2 w_i(t)^2 + \alpha a_i^2 (Cx(t))^2 \right\} \\
& \leq \frac{1}{\alpha} \sum_{i \geq N+1} \lambda_{i}^2 w_i(t)^2 + \alpha \Vert \mathcal{P}_N a \Vert_{L^2}^2 x(t)^\top C^\top C x(t) .
\end{align*}
Similarly, we obtain that
\begin{align*}
& 2 \sum_{i \geq N+1} \lambda_i w_i(t) b_i C ( A + \mu_m BC ) x(t) \leq \frac{1}{\alpha} \sum_{i \geq N+1} \lambda_{i}^2 w_i(t)^2  \\
& + \alpha \Vert \mathcal{P}_N b \Vert_{L^2}^2 x(t)^\top ( A + \mu_m BC )^\top C^\top C ( A + \mu_m BC ) x(t) ,
\end{align*}
\begin{align*}
& 2 \sum_{i \geq N+1} \lambda_i w_i(t) b_i CBC_N W(t) \leq \frac{1}{\alpha} \sum_{i \geq N+1} \lambda_{i}^2 w_i(t)^2  \\
& + \alpha \Vert \mathcal{P}_N b \Vert_{L^2}^2 W(t)^\top C_N^\top B^\top C^\top C B C_N W(t) ,
\end{align*}
and
\begin{align*}
& 2 \sum_{i \geq N+1} \lambda_i w_i(t) b_i CB \mathcal{R}(t) \\
& \leq \frac{1}{\alpha} \sum_{i \geq N+1} \lambda_{i}^2 w_i(t)^2 + \alpha \Vert \mathcal{P}_N b \Vert_{L^2}^2 (CB)^2 \mathcal{R}(t)^2 .
\end{align*}
The use of the four latter estimates implies that
\begin{align}
& \dot{V}(t) \leq \nonumber \\
& \begin{bmatrix} X(t) \\ \mathcal{R}(t) \end{bmatrix}^\top
\begin{bmatrix}
F^\top P + PF + \alpha H & PG \\ G^\top P & \alpha \Vert \mathcal{P}_N b \Vert_{L^{2}}^2 (CB)^2
\end{bmatrix}
\begin{bmatrix} X(t) \\ \mathcal{R}(t) \end{bmatrix} \nonumber \\
& + 2 \sum_{i \geq N+1} \lambda_i \left(-\lambda_i+q_c+\frac{2\lambda_i}{\alpha}\right) w_i(t)^2 \label{eq: time derivative V}
\end{align}
for $t > 0$. Since $\mathcal{R}(t) = \sum_{i \geq N+1} c_i w_i(t)$ with $c_i = \phi_i({\zeta_m})$, we infer that $\mathcal{R}(t)^2 \leq M_{1,\phi} \sum_{i \geq N+1} \lambda_i w_i(t)^2$. This implies for $t > 0$ that
\begin{align}
& \dot{V}(t) + 2 \eta V(t) \label{eq: dot_V intermediate computation} \\
& \leq \begin{bmatrix} X(t) \\ \mathcal{R}(t) \end{bmatrix}^\top
\Theta_{1,\eta}
\begin{bmatrix} X(t) \\ \mathcal{R}(t) \end{bmatrix}
+ 2\sum_{i \geq N+1} \lambda_i \Gamma_i w_i(t)^2 \nonumber
\end{align}
where $\Gamma_i = -\left(1-\frac{2}{\alpha}\right)\lambda_i + q_c + \eta + \frac{\beta M_{1,\phi}}{2}$. Now, since $\alpha > 2$, we have $\Gamma_i \leq \Gamma_{N+1}$ for all $i \geq N+1$. Moreover, combining $\Theta_{2,\eta} \preceq 0$ and the Schur complement, we infer that $\Gamma_{N+1} \leq 0$. Using also $\Theta_{1,\eta} \preceq 0$, we obtain that $\dot{V}(t) + 2 \eta V(t) \leq 0$ for all $t>0$. Since $\mathcal{A}^{1/2}w \in \mathcal{C}^0([0,\infty);L^2(0,1))$, the mapping $t \mapsto V(t)$ is continuous for $t \geq 0$, implying that $V(t) \leq e^{-2\eta t} V(0)$ for all $t \geq 0$. We now note from (\ref{eq: Lyapunov functionnal}) that $V(0) \leq \lambda_M(P) \Vert X(0) \Vert^2 + \sum_{i \geq N+1} \lambda_i \left< w_0 , \phi_i \right>^2$. Noting that $\Vert X(0) \Vert^2 \leq \Vert x_0 \Vert^2 + \Vert w_0 \Vert_{L^2}^2$ and using (\ref{eq: inner product Af and f}), we infer the existence of a constant $M_1 > 0$ such that $V(0) \leq M_1 \left( \Vert x_0 \Vert^2 + \Vert w_0 \Vert_{H^1}^2 \right)$. Using now (\ref{eq: inner product Af and f}) and (\ref{eq: Lyapunov functionnal}), we have the existence of a constant $M_2 > 0$ such that $\Vert w(t,\cdot) \Vert_{H^1}^2 \leq M_2 V(t)$. Hence, we infer the existence of a constant $M_3 > 0$ such that $\Vert w(t,\cdot) \Vert_{H^1}^2 + \Vert x(t) \Vert^2 \leq M_3 e^{-2\eta t} \left( \Vert w_0 \Vert_{H^1}^2 + \Vert x_0 \Vert^2 \right)$. The claimed conclusion follows from the change of variable (\ref{eq: change of variable 1}) and the continuous embedding $H^1(0,1) \subset L^\infty([0,1])$.
\qed

From the above proof we deduce the following corollary.

\begin{cor}\label{cor1}
In the context of Theorem~\ref{thm1}, the decay rate $\eta > 0$ of the stability estimate (\ref{eq: stab property 1}) is guaranteed provided the LMI conditions  $\Theta_{1,\eta} \preceq 0$ and $\Theta_{2,\eta} \preceq 0$ are feasible.
\end{cor}

\begin{rem}
For a given order $N$, the implementation of the conditions $\Theta_1 \prec 0$ and $\Theta_2 \prec 0$ from Theorem~\ref{thm1} require the computation of the eigenstructures $\lambda_n$ and $\phi_n$ for $1 \leq n \leq N$ as well as (an upper estimate of) $M_{1,\phi} = \sum_{i \geq N+1} \frac{\phi_i(\zeta_m)^2}{\lambda_i}$. In the case that the eigenstructures cannot be computed analytically, numerical methods can be used to estimate the $N$ first eigenstructures. Moreover, an upper bound of $M_{1,\phi} = \sum_{i \geq N+1} \frac{\phi_i(\zeta_m)^2}{\lambda_i}$ can be obtained using (\ref{eq: estimation lambda_n}) and by computing an upper bound of $\sup_{n \geq N+1} \max_{x\in[0,1]} \phi_n(x)$ by proceeding as in \cite{orlov2017general}.
\end{rem}

\subsection{Numerical illustration}

We illustrate the results of Theorem~\ref{thm1} and Corollary~\ref{cor1} for the coupled PDE-ODE system described by (\ref{eq: coupled PDE-ODE system 1}) with $\theta_1 = \pi/2$, $\theta_2 = 0$, $p=1$, $\tilde{q}=-3$, $\zeta_m = 1/4$,
\begin{align*}
A & = \begin{bmatrix}
0 & -1/4 & -1/5 & 1/5 & 1/6 \\
1/2 & 1 & -4 & 9/2 & 7/2 \\
-9/4 & -1/2 & -14 & 23 & 16 \\
-1/5 & -1/2 & -11/4 & 1/10 & 5/4 \\
-4/3 & -4/3 & -9 & 9 & 5/2
\end{bmatrix} , \\
B & = \begin{bmatrix}
-7/2 & -3/2 & -1/10 & 1/2 & 1
\end{bmatrix}^\top , \\
C & = \begin{bmatrix}
-1/10 & -1/3 & -4 & 7/8 & 7/8 
\end{bmatrix} .
\end{align*}
In this case, both PDE and ODE systems are open-loop unstable. Indeed, the dominant eigenvalue of the PDE is located approximately at $+0.533$ while the matrix $A$ has two unstable eigenvalues located approximately at $+1.046$ and $+0.247$.

\begin{figure}
     \centering
     	\subfigure[State of the reaction-diffusion system $z(t,\xi)$]{
		\includegraphics[width=3.5in]{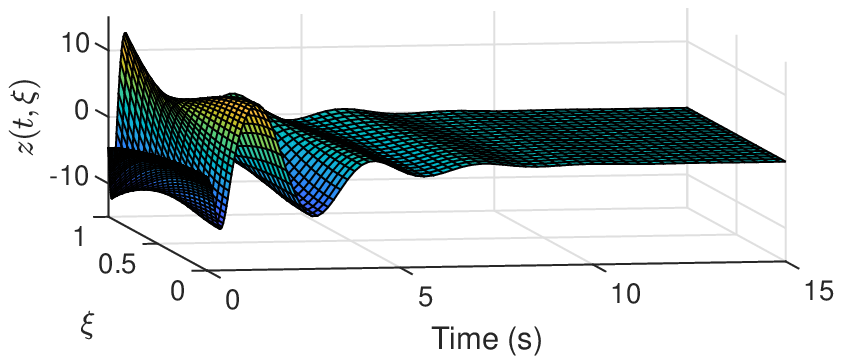}
		}
     	\subfigure[State of the ODE $x(t)$]{
		\includegraphics[width=3.5in]{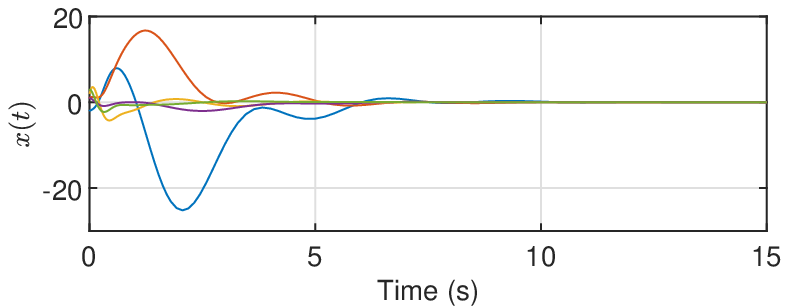}
		}
     	\subfigure[Coupling channels]{
		\includegraphics[width=3.5in]{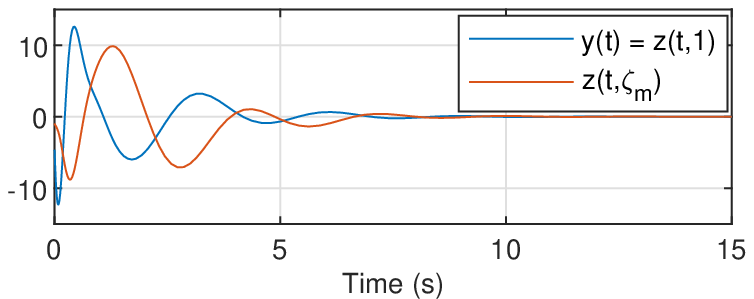}
		}		
     \caption{Time evolution of the coupled PDE-ODE system (\ref{eq: coupled PDE-ODE system 1})}
     \label{fig: sim Dirichlet actuation and measurement}
\end{figure}

We select $q = 1$ and $q_c = 4$ which satisfy (\ref{eq: decomposition of the reaction term}). Hence, we obtain that $\lambda_n = p(n-1/2)^2\pi^2+q$ and $\phi_n(\xi) = \sqrt{2}\cos((n-1/2)\pi\xi)$. Using the integral test for convergence, we infer that $M_{1,\phi} \leq \frac{2}{p \pi^2 (N-1/2)}$. The application of Theorem~\ref{thm1} with $N=3$ shows the exponential stability of the coupled PDE-ODE system (\ref{eq: coupled PDE-ODE system 1}). Moreover, the application of Corollary~\ref{cor1} with $N = 9$ shows the exponential stability of the coupled PDE-ODE system with decay rate $\eta = 0.5$. We illustrate this result with a numerical simulation. The numerical scheme consists in the modal approximation of the PDE plant by its 100 dominant modes. The initial condition is set as $w_0(\xi) = -1 + \xi^2$ and $x_0 = \begin{bmatrix} -2 & 1 & 2 & 1 & 3 \end{bmatrix}^\top$. The obtained results are depicted on Fig.~\ref{fig: sim Dirichlet actuation and measurement}, confirming the theoretical predictions of Theorem~\ref{thm1} and Corollary~\ref{cor1}.

\section{Neumann trace as an input of the ODE}\label{sec: Neumann trace as an input of the ODE}

\subsection{Coupled PDE-ODE systems}

We consider in this section the case of a reaction-diffusion PDE entering into the ODE by means of a Neumann trace instead of a Dirichlet trace. 
\begin{subequations}\label{eq: coupled PDE-ODE system 3}
\begin{align}
& z_t(t,\xi) = (p z_\xi)_\xi(t,\xi) - \tilde{q}(\xi) z(t,\xi) \\
& \cos(\theta_1) z(t,0) - \sin(\theta_1) z_\xi(t,0) = 0 \label{eq: coupled PDE-ODE system 3 - BC1} \\
& \cos(\theta_2) z(t,1) + \sin(\theta_2) z_\xi(t,1) = y(t) = C x(t) \label{eq: coupled PDE-ODE system 3 - BC2} \\
& \dot{x}(t) = A x(t) + Bz_\xi(t,\zeta_m) \label{eq: coupled PDE-ODE system 3 - ODE} \\
& z(0,\xi) = z_0(\xi) , \quad x(0) = x_0
\end{align}
\end{subequations}
for $t > 0$ and $\xi\in(0,1)$ where $\theta_1,\theta_2 \in [0,\pi/2]$, $p \in \mathcal{C}^2([0,1])$ with $p > 0$, $\tilde{q} \in \mathcal{C}^0([0,1])$, and $\zeta_m \in[0,1]$. Here $A\in\R^{n \times n}$, $B\in\R^{n}$, and $C\in\R^{1 \times n}$ are matrices, $z_0 \in L^2(0,1)$ and $x_0 \in\R^n$ are initial conditions, and $z(t,\cdot)\in L^2(0,1)$ and $x(t) \in\R^n$ are the state of the reaction-diffusion PDE and of the ODE at time $t$, respectively. 

Comparing to the PDE-ODE system (\ref{eq: coupled PDE-ODE system 1}) studied in the previous section, the PDE-ODE system (\ref{eq: coupled PDE-ODE system 3}) differs by the fact that the input of the ODE (\ref{eq: coupled PDE-ODE system 3 - ODE}) is now the pointwise Neumann trace $z_\xi(t,\zeta_m)$. In this context, the objective of this section is also to derive sufficient conditions ensuring the exponential stability of the PDE-ODE system (\ref{eq: coupled PDE-ODE system 3}) when evaluating the PDE trajectory in $H^1$-norm.

As in the previous section, we introduce without loss of generality a function $q \in \mathcal{C}^0([0,1])$ and a constant $q_c \in \R$ such that (\ref{eq: decomposition of the reaction term}) holds.

\subsection{Preliminary spectral reduction}

Considering the change of variable (\ref{eq: change of variable 1}), we infer that (\ref{eq: coupled PDE-ODE system 3}) is equivalent to
\begin{subequations}\label{eq: homogeneous coupled PDE-ODE system 3}
\begin{align}
& w_t(t,\xi) = (p w_\xi)_\xi(t,\xi) + (q_c-q(\xi)) w(t,\xi) \\
& \phantom{w_t(t,\xi) =}\; + a(\xi) y(t) + b(\xi) \dot{y}(t) \nonumber \\
& \cos(\theta_1) w(t,0) - \sin(\theta_1) w_\xi(t,0) = 0 \\
& \cos(\theta_2) w(t,1) + \sin(\theta_2) w_\xi(t,1) = 0 \\
& y(t)=Cx(t) \\
& \dot{x}(t) = A x(t) + B \left( w_\xi(t,\zeta_m) +\mu_m y(t) \right) \\
& w(0,\xi) = w_0(\xi) , \quad x(0) = x_0
\end{align}
\end{subequations}
where $a$, $b$, and $w_0$ are defined as in the previous section while $\mu_m = - b'(\zeta_m)$. Note that, after this change of variable, the well-posedness in terms of classical solutions of the above PDE-ODE systems for initial conditions $w_0 \in \cup_{\alpha_0 \in (3/4,1)} D(\mathcal{A}^{\alpha_0})$ and $x_0 \in\R^n$ is a consequence of \cite[Thm.~6.3.1 and~6.3.3]{pazy2012semigroups}. More precisely, for a given $\alpha_0 \in (3/4,1)$, we have for any $w_0 \in D(\mathcal{A}^{\alpha_0})$ and any $x_0 \in\R^n$ the existence and uniqueness of a classical solution $(w,x) \in \mathcal{C}^0([0,\infty);L^2(0,1) \times \R^n) \cap \mathcal{C}^1((0,\infty);L^2(0,1) \times \R^n)$ with $w(t,\cdot) \in D(\mathcal{A})$ for all $t > 0$. Moreover, from the proof of ~\cite[Thm.~6.3.1]{pazy2012semigroups}, we have $\mathcal{A} w \in \mathcal{C}^0((0,\infty);L^2(0,1))$ and $\mathcal{A}^{\alpha_0}w \in \mathcal{C}^0([0,\infty);L^2(0,1))$ hence $\mathcal{A}^{1/2}w \in \mathcal{C}^0([0,\infty);L^2(0,1))$.

Proceeding now as in the previous section while replacing the definition of $c_i$ by $c_i = \phi'_i(\zeta_m)$ for all $i \geq 1$, we infer that the truncated model (\ref{eq: truncated model 1}) holds while the residual dynamics is described by (\ref{eq: residual dynamics 1}).

We finally define for any $\epsilon\in(0,1/2]$ the constant $M_{2,\phi}(\epsilon) = \sum_{i \geq N+1} \frac{\phi_i'(\zeta_m)^2}{\lambda_i^{3/2+\epsilon}}$ which is finite because (\ref{eq: estimation lambda_n}) along with $\phi'_n (\zeta_m) = O(\sqrt{\lambda_n})$ as $n \rightarrow + \infty$.

\subsection{Main result}

We can now introduce the main result of this section.

\begin{thm}\label{thm2}
Let $\theta_1,\theta_2 \in [0,\pi/2]$, $p \in \mathcal{C}^2([0,1])$ with $p > 0$, $\tilde{q} \in \mathcal{C}^0([0,1])$, $\zeta_m \in[0,1]$, $A\in\R^{n \times n}$, $B\in\R^{n}$, and $C\in\R^{1 \times n}$ be given. Let $q \in \mathcal{C}^0([0,1])$ and $q_c \in\R$ be such that (\ref{eq: decomposition of the reaction term}) holds. Assume that there exist $N \geq 1$, $\epsilon\in(0,1/2]$, $P \succ 0$, $\alpha > 2$, and $\beta > 0$ such that $\Theta_1 \prec 0$, $\Theta_2 \prec 0$, and $\Theta_3 \succ 0$ where
\begin{align*}
\Theta_1 & = 
\begin{bmatrix} F^\top P + P F + \alpha H & P G \\ G^\top P & \alpha \Vert \mathcal{P}_N b \Vert_{L^2}^2 (CB)^2 - \beta \end{bmatrix} , \\
\Theta_2 & = 
\begin{bmatrix} -\lambda_{N+1}+q_c+\frac{\beta M_{2,\phi}(\epsilon)}{2}\lambda_{N+1}^{1/2+\epsilon} & \sqrt{2\lambda_{N+1}} \\ \sqrt{2\lambda_{N+1}} & -\alpha \end{bmatrix} , \\
\Theta_3 & = 
\begin{bmatrix} 1-\frac{\beta M_{2,\phi}(\epsilon)}{2\lambda_{N+1}^{1/2-\epsilon}} & \sqrt{2} \\ \sqrt{2} & \alpha \end{bmatrix} .
\end{align*} 
Then there exist constants $\eta,M > 0$ such that, for any initial conditions $z_0 \in H^2(0,1)$ and $x_0 \in \R^n$ such that $\cos(\theta_1) z_0(0) - \sin(\theta_1) z_0'(0) = 0$ and $\cos(\theta_2) z_0(1) + \sin(\theta_2) z_0'(1)=Cx_0$, the classical solution of (\ref{eq: coupled PDE-ODE system 3}) satisfies 
\begin{subequations}\label{eq: stab property 2}
\begin{equation}\label{eq: stab property 2a}
\Vert z(t,\cdot) \Vert_{H^1}^2 + \Vert x(t) \Vert^2 \leq M e^{-2\eta t} ( \Vert z_0 \Vert_{H^1}^2 + \Vert x_0 \Vert^2 )
\end{equation}
with coupling channels such that
\begin{equation}\label{eq: stab property 2b}
z_\xi(t,\zeta_m)^2 + y(t)^2 \leq  M e^{-2\eta t} ( \Vert z_0 \Vert_{H^1}^2 + \Vert \mathcal{A} w_0 \Vert_{L^2}^2 + \Vert x_0 \Vert^2 )
\end{equation}
\end{subequations}
for all $t > 0$.
\end{thm}

\begin{rem}
For any fixed $\alpha_0 \in (3/4,1)$, the estimate (\ref{eq: stab property 2a}) of Theorem~\ref{thm2} actually hold for any initial conditions such that $w_0 \in D(\mathcal{A}^{\alpha_0})$. In this case, the estimate regarding the coupling channels also holds when replacing $\Vert \mathcal{A} w_0 \Vert_{L^2}^2$ by $\Vert \mathcal{A}^{\alpha_0} w_0 \Vert_{L^2}^2$.
\end{rem}

\textbf{Proof.}
Let $N \geq 1$, $\epsilon\in(0,1/2]$, $P \succ 0$, $\alpha > 2$, and $\beta > 0$ such that $\Theta_1 \prec 0$, $\Theta_2 \prec 0$, and $\Theta_3 \succ 0$. Hence, there exist $\eta > 0$ such that $\Theta_{1,\eta} \preceq 0$ and $\Theta_{2,\eta} \preceq 0$ where
\begin{align*}
\Theta_{1,\eta} & = 
\begin{bmatrix} F^\top P + P F + 2 \eta P + \alpha H & P G \\ G^\top P & \alpha \Vert \mathcal{P}_N b \Vert_{L^2}^2 (CB)^2 - \beta \end{bmatrix} , \\
\Theta_{2,\eta} & = 
\begin{bmatrix} -\lambda_{N+1}+q_c+\eta+\frac{\beta M_{2,\phi}(\epsilon)}{2}\lambda_{N+1}^{1/2+\epsilon} & \sqrt{2\lambda_{N+1}} \\ \sqrt{2\lambda_{N+1}} & -\alpha \end{bmatrix} .
\end{align*} 
Considering the Lyapunov functionnal candidate (\ref{eq: Lyapunov functionnal}) with $X\in\R^{N+n}$ and $w \in D(\mathcal{A})$ and adopting the same approach as the one reported in the previous section, the computation of the time derivative along the system trajectories (\ref{eq: truncated model 1}) and (\ref{eq: residual dynamics 1}) gives (\ref{eq: time derivative V}) for all $t > 0$. Since $\mathcal{R}(t) = \sum_{i \geq N+1} c_i w_i(t)$ with $c_i = \phi_i'(\zeta_m)$, we infer that $\mathcal{R}(t)^2 \leq M_{2,\phi}(\epsilon) \sum_{i \geq N+1} \lambda_i^{3/2+\epsilon} w_i(t)^2$. This implies that (\ref{eq: dot_V intermediate computation}) holds for $t > 0$ with $\Gamma_i = -\left(1-\frac{2}{\alpha}\right)\lambda_i + q_c + \eta + \frac{\beta M_{2,\phi}(\epsilon)}{2}\lambda_{i}^{1/2+\epsilon}$. Since $\epsilon \in(0,1/2]$, we observe for $i \geq N+1$ that $\lambda_{i}^{1/2+\epsilon} = \lambda_i/\lambda_i^{1/2-\epsilon} \leq \lambda_i/\lambda_{N+1}^{1/2-\epsilon}$ hence $\Gamma_i \leq -\left(1-\frac{2}{\alpha}-\frac{\beta M_{2,\phi}(\epsilon)}{2\lambda_{N+1}^{1/2-\epsilon}}\right)\lambda_i + q_c + \eta$. Using $\Theta_3 \succ 0$ and Schur's complement, we infer that $\Gamma_i \leq - \lambda_{N+1} + q_c + \eta + \frac{\beta M_{2,\phi}(\epsilon)}{2}\lambda_{N+1}^{1/2+\epsilon} + \frac{2\lambda_{N+1}}{\alpha}$ for all $i \geq N+1$. Using now $\Theta_{2,\eta} \preceq 0$ and Schur's complement, we obtain that $\Gamma_i \leq 0$ for all $i \geq N+1$. Combining this result with $\Theta_{1,\eta} \preceq 0$, we deduce from (\ref{eq: dot_V intermediate computation}) that $\dot{V}(t) + 2 \eta V(t) \leq 0$ for all $t>0$. From now on, the proof of (\ref{eq: stab property 2a}) follows from the same arguments than the ones reported in the previous section. To complete the proof, we only need to establish the exponential decrease of the term $z_\xi(t,\zeta_m)$ to obtain (\ref{eq: stab property 2b}). This is done in Apprendix by invoking a small gain argument.
\qed

\begin{cor}\label{cor2}
In the context of Theorem~\ref{thm2}, the decay rate $\eta > 0$ of the stability estimate (\ref{eq: stab property 2}) is guaranteed provided the LMI conditions  $\Theta_{1,\eta} \preceq 0$, $\Theta_{2,\eta} \preceq 0$, and $\Theta_3 \succ 0$ are feasible.
\end{cor}

\subsection{Numerical illustration}

We illustrate the results of Theorem~\ref{thm2} and Corollary~\ref{cor2} for the coupled PDE-ODE system described by (\ref{eq: coupled PDE-ODE system 3}) with $\theta_1 = 0$, $\theta_2 = \pi/2$, $p=1$, $\tilde{q}=-3$, $\zeta_m = 1/4$,
\begin{align*}
A & = \begin{bmatrix}
-1/4 & -1/6 & 2 & 1 & 1/12 \\
-3/2 & -3/2 & 5 & 5 & 1/6 \\
3/2 & -4 & -15/2 & -5 & -1/3 \\
-13/2 & 22 & 22 & -14 & -1/2 \\
1/7 & -1/2 & -1/2 & 1/5 & -5/2 
\end{bmatrix} , \\
B & = \begin{bmatrix}
-5/4 & 2/3 & 1/6 & -1/6 & 0
\end{bmatrix}^\top , \\
C & = \begin{bmatrix}
-2/5 & -5/4 & 3/2 & 1/3 & 1/40 
\end{bmatrix} .
\end{align*}
Both PDE and ODE systems are open-loop unstable. Indeed, the dominant eigenvalue of the PDE is located approximately at $+0.533$ while the matrix $A$ has one unstable eigenvalue located approximately at $+0.393$. 

\begin{figure}
     \centering
     	\subfigure[State of the reaction-diffusion system $z(t,\xi)$]{
		\includegraphics[width=3.5in]{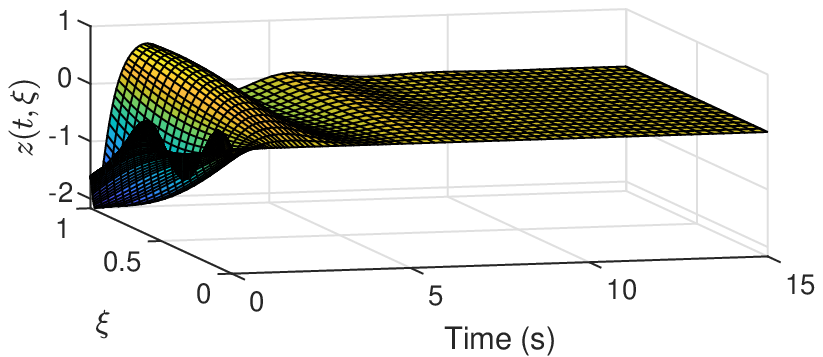}
		}
     	\subfigure[State of the ODE $x(t)$]{
		\includegraphics[width=3.5in]{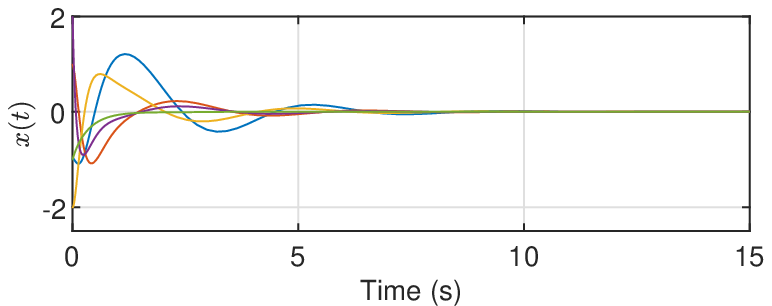}
		}
     	\subfigure[Coupling channels]{
		\includegraphics[width=3.5in]{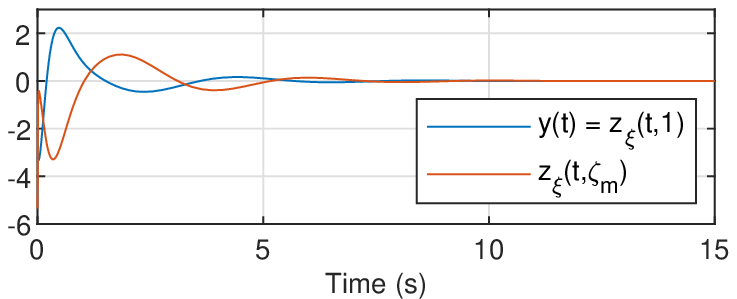}
		}		
     \caption{Time evolution of the coupled PDE-ODE system (\ref{eq: coupled PDE-ODE system 3})}
     \label{fig: sim Neumann actuation and measurement}
\end{figure}

We select $q = 1$ and $q_c = 4$ which satisfy (\ref{eq: decomposition of the reaction term}). Hence, we obtain that $\lambda_n = p(n-1/2)^2\pi^2+q$ and $\phi_n(\xi) = \sqrt{2}\sin((n-1/2)\pi\xi)$. Using the integral test for convergence, we infer that $M_{2,\phi}(\epsilon) \leq \frac{1}{\epsilon p^{3/2+\epsilon} \pi^{1+2\epsilon} (N-1/2)^{2\epsilon}}$. The application of Theorem~\ref{thm2} with $\epsilon = 1/6$ and $N=2$ shows the exponential stability of the coupled PDE-ODE system (\ref{eq: coupled PDE-ODE system 1}). Moreover, the application of Corollary~\ref{cor2} with $N = 10$ shows the exponential stability of the coupled PDE-ODE system with decay rate $\eta = 0.4$. We illustrate this result with a numerical simulation. The numerical scheme consists in the modal approximation of the PDE plant by its 100 dominant modes. The initial condition is set as
$w_0(\xi) = 5\xi(1-\xi)^2\cos(3\pi\xi)$ and $x_0 = \begin{bmatrix} -1 & 1 & -2 & 2 & -1 \end{bmatrix}^\top$. The obtained results are depicted on Fig.~\ref{fig: sim Neumann actuation and measurement}, confirming the theoretical predictions of Theorem~\ref{thm2} and Corollary~\ref{cor2}.

\section{Conclusion}\label{sec: conclusion}
This paper has addressed the topic of assessing the stability of coupled systems composed of a reaction-diffusion equation and a finite-dimensional linear time-invariant ODE. The considered coupling channels are located either at the boundaries or in the domain of the PDE and consist of the input and output signals of the ODE. The reported sufficient stability conditions take the form of tractable LMIs and have been derived by adopting a spectral reduction-based method. Moreover, we have also assessed the exponential decrease to zero of the aforementioned coupling channels, particularly in the case of Neumann boundary couplings. The drawback of the present Lyapunov function based approach is that the derived stability condition are only sufficient, hence may be not satisfied by some stable reaction-diffusion systems. Nevertheless, as illustrated via the reported numerical examples, this method can be successfully applied to assess the exponential stability of coupled PDE-ODE systems for which both the open-loop PDE and ODE plants are exponentially unstable.


\bibliographystyle{plain}        
\bibliography{autosam}           

\begin{thebibliography}{10}

\bibitem{ahmadi2016dissipation}
Mohamadreza Ahmadi, Giorgio Valmorbida, and Antonis Papachristodoulou.
\newblock Dissipation inequalities for the analysis of a class of {PDEs}.
\newblock {\em Automatica}, 66:163--171, 2016.

\bibitem{ayadi2018exponential}
Habib Ayadi.
\newblock Exponential stabilization of cascade {ODE-linearized KdV} system by
  boundary dirichlet actuation.
\newblock {\em European Journal of Control}, 43:33--38, 2018.

\bibitem{barreau2018input}
Matthieu Barreau, Fr{\'e}d{\'e}ric Gouaisbaut, Alexandre Seuret, and Rifat
  Sipahi.
\newblock Input/output stability of a damped string equation coupled with
  ordinary differential system.
\newblock {\em International Journal of Robust and Nonlinear Control},
  28(18):6053--6069, 2018.

\bibitem{matthieu2020integral}
Matthieu Barreau, Carsten Scherer, Fr{\'e}d{\'e}ric Gouaisbaut, and Alexandre
  Seuret.
\newblock Integral quadratic constraints on linear infinite-dimensional systems
  for robust stability analysis.
\newblock In {\em IFAC World Congress}, 2020.

\bibitem{barreau2018exponential}
Matthieu Barreau, Alexandre Seuret, and Fr{\'e}d{\'e}ric Gouaisbaut.
\newblock Exponential {Lyapunov} stability analysis of a drilling mechanism.
\newblock In {\em 2018 IEEE Conference on Decision and Control (CDC)}, pages
  6579--6584. IEEE, 2018.

\bibitem{barreau2018lyapunov}
Matthieu Barreau, Alexandre Seuret, Fr{\'e}d{\'e}ric Gouaisbaut, and Lucie
  Baudouin.
\newblock Lyapunov stability analysis of a string equation coupled with an
  ordinary differential system.
\newblock {\em IEEE Transactions on Automatic Control}, 63(11):3850--3857,
  2018.

\bibitem{baudouin2019stability}
Lucie Baudouin, Alexandre Seuret, and Fr{\'e}d{\'e}ric Gouaisbaut.
\newblock Stability analysis of a system coupled to a heat equation.
\newblock {\em Automatica}, 99:195--202, 2019.

\bibitem{bekiaris2010compensating}
Nikolaos Bekiaris-Liberis and Miroslav Krstic.
\newblock Compensating the distributed effect of a wave pde in the actuation or
  sensing path of {MIMO LTI} systems.
\newblock {\em Systems \& Control Letters}, 59(11):713--719, 2010.

\bibitem{boulouz2021well}
Abed Boulouz, Hamid Bounit, and Said Hadd.
\newblock Well-posedness and exponential stability of boundary control systems
  with dynamic boundary conditions.
\newblock {\em Systems \& Control Letters}, 147:104825, 2021.

\bibitem{coron2004global}
Jean-Michel Coron and Emmanuel Tr{\'e}lat.
\newblock Global steady-state controllability of one-dimensional semilinear
  heat equations.
\newblock {\em SIAM Journal on Control and Optimization}, 43(2):549--569, 2004.

\bibitem{coron2006global}
Jean-Michel Coron and Emmanuel Tr{\'e}lat.
\newblock Global steady-state stabilization and controllability of {1D}
  semilinear wave equations.
\newblock {\em Commun. Contemp. Math.}, 8(04):535--567, 2006.

\bibitem{daafouz2014nonlinear}
Jamal Daafouz, Marius Tucsnak, and Julie Valein.
\newblock Nonlinear control of a coupled {PDE/ODE} system modeling a switched
  power converter with a transmission line.
\newblock {\em Systems \& Control Letters}, 70:92--99, 2014.

\bibitem{das2020robust}
Amritam Das, Sachin Shivakumar, Matthew Peet, and Siep Weiland.
\newblock Robust analysis of uncertain {ODE-PDE} systems using {PI}
  multipliers, {PIEs} and {LPIs}.
\newblock In {\em 2020 59th IEEE Conference on Decision and Control (CDC)},
  pages 634--639. IEEE, 2020.

\bibitem{gahlawat2016convex}
Aditya Gahlawat and Matthew~M Peet.
\newblock A convex sum-of-squares approach to analysis, state feedback and
  output feedback control of parabolic {PDEs}.
\newblock {\em IEEE Transactions on Automatic Control}, 62(4):1636--1651, 2016.

\bibitem{he2016cooperative}
Wei He and Shuzhi~Sam Ge.
\newblock Cooperative control of a nonuniform gantry crane with constrained
  tension.
\newblock {\em Automatica}, 66:146--154, 2016.

\bibitem{krstic2009compensatingString}
Miroslav Krstic.
\newblock Compensating a string {PDE} in the actuation or sensing path of an
  unstable {ODE}.
\newblock {\em IEEE Transactions on Automatic Control}, 54(6):1362--1368, 2009.

\bibitem{krstic2009compensating}
Miroslav Krstic.
\newblock Compensating actuator and sensor dynamics governed by diffusion
  {PDEs}.
\newblock {\em Systems \& Control Letters}, 58(5):372--377, 2009.

\bibitem{krstic2010compensation}
Miroslav Krstic and Nikolaos Bekiaris-Liberis.
\newblock Compensation of infinite-dimensional input dynamics.
\newblock {\em Annual Reviews in Control}, 34(2):233--244, 2010.

\bibitem{lhachemi2020feedback}
Hugo Lhachemi and Christophe Prieur.
\newblock Feedback stabilization of a class of diagonal infinite-dimensional
  systems with delay boundary control.
\newblock {\em IEEE Transactions on Automatic Control}, 66(1):105--120, 2020.

\bibitem{lhachemi2018boundary}
Hugo Lhachemi, David Saussi{\'e}, and Guchuan Zhu.
\newblock Boundary feedback stabilization of a flexible wing model under
  unsteady aerodynamic loads.
\newblock {\em Automatica}, 97:73--81, 2018.

\bibitem{li2014adaptive}
Jian Li and Yungang Liu.
\newblock Adaptive stabilization of coupled {PDE--ODE} systems with multiple
  uncertainties.
\newblock {\em ESAIM: Control, Optimisation and Calculus of Variations},
  20(2):488--516, 2014.

\bibitem{natarajan2021compensating}
Vivek Natarajan.
\newblock Compensating {PDE} actuator and sensor dynamics using {Sylvester}
  equation.
\newblock {\em Automatica}, 123:109362, 2021.

\bibitem{nicaise2009stability}
Serge Nicaise, Julie Valein, and Emilia Fridman.
\newblock Stability of the heat and of the wave equations with boundary
  time-varying delays.
\newblock {\em Discrete and Continuous Dynamical Systems?` Series S}, 2(3):559,
  2009.

\bibitem{orlov2017general}
Yury Orlov.
\newblock On general properties of eigenvalues and eigenfunctions of a
  {Sturm--Liouville} operator: comments on ''{ISS} with respect to boundary
  disturbances for {1-D} parabolic {PDEs}''.
\newblock {\em IEEE Transactions on Automatic Control}, 62(11):5970--5973,
  2017.

\bibitem{pazy2012semigroups}
Amnon Pazy.
\newblock {\em Semigroups of linear operators and applications to partial
  differential equations}, volume~44.
\newblock Springer Science \& Business Media, 2012.

\bibitem{peet2021partial}
Matthew~M Peet.
\newblock A {Partial Integral Equation (PIE)} representation of coupled linear
  {PDEs} and scalable stability analysis using {LMIs}.
\newblock {\em Automatica}, 125:109473, 2021.

\bibitem{peet2021representation}
Matthew~M Peet.
\newblock Representation of networks and systems with delay: {DDEs}, {DDFs},
  {ODE--PDEs} and {PIEs}.
\newblock {\em Automatica}, 127:109508, 2021.

\bibitem{sanz2018robust}
Ricardo Sanz, Pedro Garc{\'\i}a, and Miroslav Krstic.
\newblock Robust compensation of delay and diffusive actuator dynamics without
  distributed feedback.
\newblock {\em IEEE Transactions on Automatic Control}, 64(9):3663--3675, 2018.

\bibitem{shivakumar2019generalized}
Sachin Shivakumar, Amritam Das, Siep Weiland, and Matthew~M Peet.
\newblock A generalized {LMI} formulation for input-output analysis of linear
  systems of {ODEs} coupled with {PDEs}.
\newblock In {\em 2019 IEEE 58th Conference on Decision and Control (CDC)},
  pages 280--285. IEEE, 2019.

\bibitem{susto2010control}
Gian~Antonio Susto and Miroslav Krstic.
\newblock Control of {PDE--ODE} cascades with neumann interconnections.
\newblock {\em Journal of the Franklin Institute}, 347(1):284--314, 2010.

\bibitem{tang2011state}
Shuxia Tang and Chengkang Xie.
\newblock State and output feedback boundary control for a coupled {PDE--ODE}
  system.
\newblock {\em Systems \& Control Letters}, 60(8):540--545, 2011.

\bibitem{wang2019output}
Ji~Wang and Miroslav Krstic.
\newblock Output feedback boundary control of a heat {PDE} sandwiched between
  two {ODEs}.
\newblock {\em IEEE Transactions on Automatic Control}, 64(11):4653--4660,
  2019.

\bibitem{wang2015sliding}
Jun-Min Wang, Jun-Jun Liu, Beibei Ren, and Jinhao Chen.
\newblock Sliding mode control to stabilization of cascaded heat {PDE--ODE}
  systems subject to boundary control matched disturbance.
\newblock {\em Automatica}, 52:23--34, 2015.

\bibitem{wu2014static}
Huai-Ning Wu and Jun-Wei Wang.
\newblock Static output feedback control via {PDE} boundary and {ODE}
  measurements in linear cascaded {ODE--beam} systems.
\newblock {\em Automatica}, 50(11):2787--2798, 2014.

\end{thebibliography}



\appendix
\section{End of the proof of Theorem~\ref{thm2}}
We investigate the exponential decrease of $z_\xi(t,\zeta_m)$ to zero. Using the change of variable (\ref{eq: change of variable 1}) and the identity $y(t) = Cx(t)$, we have $z_\xi(t,\zeta_m) = w_\xi(t,\zeta_m) + \mu_m Cx(t)$ for $t > 0$. Hence, based on (\ref{eq: stab property 2a}), we only need to study the term $w_\xi(t,\zeta_m) = \sum_{j \geq 1} \phi_j'(\zeta_m) w_j(t)$. Let $\alpha_0 \in (3/4,1)$ and $w_0 \in D(\mathcal{A}^{\alpha_0})$. Let $N_0 \geq 1$ and $\kappa > 0$ be such that $-\lambda_n + q_c \leq - \eta - \kappa$ for all $n \geq N_0 + 1$. Consider an arbitrary fixed integer $M \geq N_0$. Then we have
\begin{align*}
\vert w_\xi(t,\zeta_m) \vert 
& \leq \sum_{j \geq 1} \vert \phi_j'(\zeta_m) \vert \vert  w_j(t) \vert \\
& \leq C_M \sqrt{\sum_{j = 1}^M w_j(t)^2} + C_\phi \sqrt{\sum_{j \geq M+1} \lambda_j^{2\alpha_0} w_j(t)^2} \\
&\leq C_M \Vert w(t,\cdot) \Vert_{L^2} + C_\phi \sqrt{\sum_{j \geq M+1} \lambda_j^{2\alpha_0} w_j(t)^2}
\end{align*}
where $C_M = \sqrt{\sum_{j = 1}^M \phi_j'(\zeta_m)^2}$ and $C_\phi = \sqrt{\sum_{j \geq 1} \frac{\phi_j'(\zeta_m)^2}{\lambda_j^{2\alpha_0}}} < \infty$. Based again on (\ref{eq: change of variable 1}) and (\ref{eq: stab property 2a}) we only need to study the term $S_M(t) = \sum_{j \geq M+1} \lambda_j^{2\alpha_0} w_j(t)^2$. To do so, we integrate (\ref{eq: spectral reduction 1 - 1}) for $i \geq N_0 +1$ and direct estimations give
\begin{align*}
& \vert \lambda_i^{\alpha_0} w_i(t) \vert 
\leq e^{(-\lambda_i+q_c)t} \vert \lambda_i^{\alpha_0} w_i(0) \vert \\
& + \left( \vert a_i \vert \Vert C \Vert + \vert b_i \vert \Vert CA_e \Vert \right) \int_0^t \lambda_i^{\alpha_0} e^{(-\lambda_i+q_c)(t-s)} \Vert x(s) \Vert \diff s \\
& + \vert b_i \vert \vert CB \vert \int_0^t \lambda_i^{\alpha_0} e^{(-\lambda_i+q_c)(t-s)} \vert w_\xi(s,\zeta_m) \vert \diff s 
\end{align*}
with $A_e = A + \mu_m BC$. For any $f \in L^\infty_\mathrm{loc}(\R)$, we have 
\begin{align*}
& \int_0^t \lambda_i^{\alpha_0} e^{(-\lambda_i+q_c)(t-s)} \vert f(s) \vert \diff s \\
& = e^{-\eta t} \int_0^t \lambda_i^{\alpha_0} e^{(-\lambda_i+q_c+\eta)(t-s)} \times e^{\eta s} \vert f(s) \vert \diff s \\
& \leq \frac{\lambda_i^{\alpha_0}}{\lambda_i-q_c-\eta} e^{-\eta t} \mathrm{ess\,sup}_{s \in [0,t]} e^{\eta s} \vert f(s) \vert 
\end{align*}
because $\lambda_i - q_c - \eta \geq \kappa > 0$ for all $i \geq N_0 +1$. Since $\alpha_0 \in (3/4,1)$ we also have that $\lambda_i^{\alpha_0}/(\lambda_i-q_c-\eta) \rightarrow 0$ as $i \rightarrow + \infty$. Hence there exists a constant $\sigma = \sigma(\alpha_0) > 0$, independent of $M$, so that $\lambda_i^{\alpha_0}/(\lambda_i-q_c-\eta) \leq \sigma$ for all $i \geq N_0+1$. Combining the latter estimates, we infer that
\begin{align*}
& \vert \lambda_i^{\alpha_0} w_i(t) \vert 
\leq e^{-\eta t} \vert \lambda_i^{\alpha_0} w_i(0) \vert \\
& + \sigma \left( \vert a_i \vert \Vert C \Vert + \vert b_i \vert \Vert CA_e \Vert \right) e^{-\eta t} \sup_{s \in [0,t]} e^{\eta s} \Vert x(s) \Vert \\
& + \sigma C_M \vert b_i \vert \vert CB \vert e^{-\eta t} \sup_{s \in [0,t]} e^{\eta s} \Vert w(s,\cdot) \Vert_{L^2} \\
& + \sigma C_\phi \vert b_i \vert \vert CB \vert e^{-\eta t} \sup_{s \in [0,t]} e^{\eta s} \sqrt{S_M(s)}  
\end{align*}
for all $i \geq N_0 +1$ and all $t \geq 0$. The use of Young's inequality and summing for $i \geq M+1 \geq N_0 +1$ we obtain that
\begin{align*}
& S_M(t) 
\leq 5 e^{- 2\eta t} S_M(0) \\
& + 5 \sigma^2 \left( \Vert \mathcal{P}_M a \Vert_{L^2}^2 \Vert C \Vert^2 + \Vert \mathcal{P}_M b \Vert_{L^2}^2 \Vert CA_e \Vert^2 \right) \\
& \phantom{+}\; \times e^{-2\eta t} \sup_{s \in [0,t]} e^{2\eta s} \Vert x(s) \Vert^2 \\
& + 5 \sigma^2 C_M^2 \Vert \mathcal{P}_M b \Vert_{L^2}^2 \vert CB \vert^2 e^{-2\eta t} \sup_{s \in [0,t]} e^{2\eta s} \Vert w(s,\cdot) \Vert_{L^2}^2 \\
& + 5 \sigma^2 C_\phi^2 \Vert \mathcal{P}_M b \Vert_{L^2}^2 \vert CB \vert^2 e^{-2\eta t} \sup_{s \in [0,t]} e^{2\eta s} S_M(s)  
\end{align*}
for all $t \geq 0$, hence 
\begin{align*}
& \sup_{s \in [0,t]} e^{2\eta s} S_M(s) 
\leq 5 S_M(0) \\
& + 5 \sigma^2 \left( \Vert \mathcal{P}_M a \Vert_{L^2}^2 \Vert C \Vert^2 + \Vert \mathcal{P}_M b \Vert_{L^2}^2 \Vert CA_e \Vert^2 \right)  \\
& \phantom{+}\; \times \sup_{s \in [0,t]} e^{2\eta s} \Vert x(s) \Vert^2 \\
& + 5 \sigma^2 C_M^2 \Vert \mathcal{P}_M b \Vert_{L^2}^2 \vert CB \vert^2 \sup_{s \in [0,t]} e^{2\eta s} \Vert w(s,\cdot) \Vert_{L^2}^2 \\
& + 5 \sigma^2 C_\phi^2 \Vert \mathcal{P}_M b \Vert_{L^2}^2 \vert CB \vert^2 \sup_{s \in [0,t]} e^{2\eta s} S_M(s) .
\end{align*}
Since $\Vert \mathcal{P}_M b \Vert_{L^2} \rightarrow 0$ when $M \rightarrow +\infty$, we infer the existence of a large enough integer $M \geq N_0 + 1$, independent of the initial condition $w_0 \in D(\mathcal{A}^{\alpha_0})$, such that $5 \sigma^2 C_\phi^2 \Vert \mathcal{P}_M b \Vert_{L^2}^2 \vert CB \vert^2 < 1$. Fixing such a $M \geq N_0 + 1$ and because all the supremums appearing in the latter estimate are finite (recall that $\mathcal{A}^{\alpha_0}w \in \mathcal{C}^0([0,\infty);L^2(0,1))$), we obtain the existence of a constant $M_4 > 0$ such that 
\begin{align*}
\sup_{s \in [0,t]} e^{2\eta s} S_M(s) 
& \leq M_4 S_M(0) + M_4 \sup_{s \in [0,t]} e^{2\eta s} \Vert x(s) \Vert^2 \\
& \phantom{\leq}\; + M_4 \sup_{s \in [0,t]} e^{2\eta s} \Vert w(s,\cdot) \Vert_{L^2}^2 
\end{align*}
for all $t \geq 0$. Noting that $S_M(0) \leq \Vert \mathcal{A}^{\alpha_0} w_0 \Vert_{L^2}^2$, the claimed conclusion follows from (\ref{eq: change of variable 1}) and (\ref{eq: stab property 2a}). In the case $w_0 \in D(\mathcal{A})$, it can easily be seen that $D(\mathcal{A}) \subset D(\mathcal{A}^{\alpha_0})$ and $\Vert \mathcal{A}^{\alpha_0} w_0 \Vert_{L^2}^2 \leq \Vert w_0 \Vert_{L^2}^2 + \Vert \mathcal{A} w_0 \Vert_{L^2}^2$, which gives (\ref{eq: stab property 2b}) and concludes the proof.

\end{document}